\documentclass[11pt]{article}
\usepackage{amssymb}
\usepackage{amsmath}
\usepackage{amsthm}
\usepackage{mathtools}
\usepackage{mathrsfs}
\usepackage{accents} 
\usepackage{etoolbox}
\usepackage{siunitx}
\usepackage{dsfont} 
\usepackage{graphicx}
\usepackage{color}
\usepackage[oldvoltagedirection]{circuitikz}
\usepackage[margin=10pt,font=small,labelfont=bf,labelsep=endash]{caption}
\usepackage{subcaption}

\numberwithin{equation}{section}

\usepackage{algorithm}
\usepackage{algorithmicx}
\usepackage{algpseudocode}

\newcommand{\node}{u}

\newcommand{\arc}{e}

\newcommand{\velocity}{v}
\newcommand{\vel}{\velocity}

\newcommand{\pressure}{p}
\newcommand{\press}{\pressure}

\newcommand{\bore}{D}
\newcommand{\diameter}{\bore}
\newcommand{\diam}{\diameter}

\newcommand{\abs}[1]{\lvert\fcarg{#1}\rvert}

\makeatletter
\newcommand{\fcdot}{\,\cdot\,}
\newcommand{\fcarg}[1]{\def\fc@rg{#1}\ifx\fc@rg\empty\fcdot\else\fc@rg\fi}
\makeatother

\newcommand{\gravity}{g}
\newcommand{\grav}{\gravity}

\newcommand{\densityinternalenergy}{e}
\newcommand{\densintene}{\densityinternalenergy}
\newcommand{\densityinternalentropy}{s}
\newcommand{\densintent}{\densityinternalentropy}

\newcommand{\heattransfercoefficient}{k_W}
\newcommand{\heattrans}{\heattransfercoefficient}

\newcommand{\temperature}{T}
\newcommand{\temp}{\temperature}

\newcommand{\soiltemp}{\temp_\text{W}}

\newcommand{\diff}{\,{d}}

\title{Mathematical modeling, simulation and {optimization via (port-Hamiltonian)} differential-algebraic equations}

\author{Volker Mehrmann\\
Institut für Mathematik, {IMoS MA 4-5}, TU Berlin\\ {Fasanenstr. 89,} 10623 Berlin, Germany\\
E-mail: {mehrmann@math.tu-berlin.de}}

\newtheorem{theorem}{Theorem}[section]

\newtheorem{definition}[theorem]{Definition}

\newtheorem{example}[theorem]{Example}
\newtheorem{assumption}[theorem]{Assumption}

\newcommand{\dd}[2]{\frac{\td#1}{\td#2}}

\newcommand{\mat}[1]{\begin{matrix}#1\end{matrix}}
\newcommand{\bmat}[1]{\begin{bmatrix}#1\end{bmatrix}}

\newcommand{\pa}{\partial}
\newcommand{\mB}[1]{\mathbb{#1}}
\newcommand{\mc}[1]{\mathcal{#1}}
\newcommand{\R}{\mB{R}}

\newcommand{\td}{\textup{d}}
\definecolor{myBlue}{RGB}{30,144,255} 
\definecolor{myGreen}{RGB}{69,169,0} 
\definecolor{myRed}{RGB}{165,12,42}
\definecolor{myOrange}{RGB}{225,92,22}
\definecolor{color0}{rgb}{0.121568,0.46666,0.70588}
\definecolor{color1}{rgb}{1,0.4980,0.0549}
\definecolor{color2}{rgb}{0.17254,0.62745,0.17254901}
\definecolor{color3}{rgb}{0.83921,0.15294,0.156862}
\definecolor{color4}{rgb}{0.580392,0.4039215,0.7411764}
\definecolor{color5}{rgb}{0,0,0}

\newcommand{\delete}[1]{ }

\def\R{\mathbb{R}}

\newcommand{\ddt}{\ensuremath{\tfrac{\mathrm{d}}{\mathrm{d}t}} }

\DeclareMathOperator{\diag}{diag}

\DeclareMathOperator{\rank}{rank}

\DeclareMathOperator{\corank}{corank}

\newcommand{\fc}[2]{\mathcal  C^{#1}(\mathcal S,#2)}

\makeatletter
\newsavebox{\@brx}
\newcommand{\llangle}[1][]{\savebox{\@brx}{\(\m@th{#1\langle}\)}%
  \mathopen{\copy\@brx\kern-0.5\wd\@brx\usebox{\@brx}}}
\newcommand{\rrangle}[1][]{\savebox{\@brx}{\(\m@th{#1\rangle}\)}%
  \mathclose{\copy\@brx\kern-0.5\wd\@brx\usebox{\@brx}}}
\makeatother

\newcommand{\calA}{\ensuremath{\mathcal{A}} }

\newcommand{\hamiltonian}{\mathcal{H}}

\begin{document}

\maketitle
\begin{abstract}
Differential-algebraic equations (DAEs) present a ubiquitous model class in the automated modeling, simulation, control and optimization of complex finite or infinite dimensional dynamical systems, in particular in the construction of digital twins. When combined with adequate space and time discretization methods they present the work-horse of developments in many engineering domains. However, this approach comes  with certain challenges for the numerical methods and in particular in the coupling of systems from different physical domains with different scales or regularity requirements. We survey the modeling with DAEs and address how to overcome some of the challenges via the framework of port-Hamiltonian differential-algebraic systems (pHDAEs) and their simulation and  optimization with space-time-model adaptive methods.
\end{abstract}

{\bf Congratulations to the Institut für Techno- und Wirtschaftsmathematik (ITWM) Kaiserslautern on the occasion of their 30th anniversary.}

\section{Introduction}
\label{sec:introduction}
In modern engineering, the  modeling, simulation, control and optimization of complex dynamical systems is a key factor in the design and operation of new products or processes. {A well established modularized and automated system theoretic approach for this is to combine a large number of submodels via a network. This is well understood in many application domains such as e.g. circuit simulation or multi-body dynamics,  and it is a driver of modern systems engineering.

A major factor in this development is the use of \emph{differential-algebraic equations (DAEs), see \cite{KunM24} for a detailed overview of the analysis, numerical simulation, and control techniques for this class. This approach is also successful for distributed parameter systems, where typically the partial differential equation model is first discretized in space to obtain a finite dimensional DAE system. In this survey we will focus on finite dimensional DAE systems.}

Furthermore, when a real product or process  is equipped with a large number of sensors that collect and store data, then one can even  use the network based approach to build a \emph{digital twin} which  uses the acquired data combined with numerical simulations to accompany  its physical twin, ideally throughout its life cycle \cite{BouVC20,FalK21,RasO20}. See \cite{DTreport24} for a detailed description of the state of the art and challenges for digital twins.

A driving force {and a real success story} in the technological development, in particular for digital twins, are accurate and efficient automated modeling and simulation tools such as e.g. SPICE\footnote{\tt{https://bwrcs.eecs.berkeley.edu/Classes/IcBook/SPICE}} in circuit simulation, SIMPACK\footnote{\tt{www.simpack.de/}} in multi-body dynamics, or Simulink\footnote{\tt{https://www.mathworks.com/products/simulink.html}} for general systems engineering.
Modeling languages like Modelica\footnote{\tt{https://doc.modelica.org/}} or  Simscape\footnote{\tt{https://www.mathworks.com/products/simscape.html}}, as well as computer algebra based packages like
MAPLE\footnote{\tt{www.maplesoft.com/}} or MATHEMATICA\footnote{\tt{www.wolfram.com/}} provide established tools to generate  dynamical system models for simulation and control.
Furthermore, for distributed parameter systems, packages for the solution of partial differential equations (PDEs) are required. Also for this task  there exist a large number commercial and open source numerical simulation packages, like e.g. Abaqus\footnote{\tt{https://www.3ds.com/products-services/simulia/products/abaqus/}}, Ansys\footnote{\tt{https://www.ansys.com/}}, COMSOL\footnote{\tt{https://www.comsol.com}}, Dymola\footnote{\tt{https://www.3ds.com/products-services/catia/products/dymola/}}, FEniCS\footnote{\tt{https://fenicsproject.org}}, and also in-house packages like e.g. AVL CRUISE${}^{TM}$M.


Despite a very successful use of network based automated modeling in many engineering domains,  the coupling of components from different physical domains and different simulation software packages still presents a major obstacle, see e.g. \cite{GomTBLV18}, due to different scalings, different software environments, as well as different accuracy and real-time requirements. This is, in particular, a challenge for the use of these tools in digital twins.

To achieve a balance between computational effort, real time requirements, and desired accuracy, one typically uses the well-known concept of adaptivity in space and time which steers the grid sizes and time steps based on error estimates \cite{BanR03,NocSV09}. {Recently, also the adaptivity in the model selection is employed, see e.g.~\cite{DomDSLM18,StoM18}. This approach uses model catalogs of coarse and fine discretizations, as well as surrogate and  reduced-order models, to switch between different  models in order to achieve the required accuracy,  see e.g. \cite{DomHLMMT21,HauMMMMRS20} in the context of energy transport networks.  This approach has been particularly successful in the context of optimization problems~\cite{DaeMRS24,KruMS21,MehSS18} to obtain feasible solutions in each step of an optimization method. Note, however that such model catalogs are typically domain specific and thus this approach still presents a challenge in the context of multi-domain models, in particular for digital twins, where the model has to be adapted in view of data from the physical twin,  see e.g.~\cite{CheSMGL26}.}

In this paper we survey the theoretical and practical aspects of  modeling, and simulation, as well as { optimization} of DAE systems, describe some of the challenges and illustrate that many of them can be overcome when combining the classical DAE approach with the approach of energy based modeling via port-Hamiltonian DAE systems (pHDAEs).

In Section~\ref{sec:generaldaes} we give a brief overview over the theory of DAE systems and discuss some of the challenges that arise in {numerical} simulation, control and {optimization}. In Section~\ref{sec:phdae} we discuss the class of port-Hamiltonian DAEs (pHDAEs) and how they can be used to overcome some of these challenges. We illustrate pHDAEs in the industrial context with simple models from circuit simulation and multi-body dynamics  and then in Section~\ref{sec:adaptive} demonstrate how they can be used in the simulation and optimization of gas transport  {and district heating networks} via space-time-model adaptivity. The paper closes with an extended summary and an outlook.

\section{General differential-algebraic equation systems}\label{sec:generaldaes}
We use the following notation. By $C^k(\mathcal D,\mathcal W)$ we denote the set of $k$-times continuously differential functions from $\mathcal D$ to $\mathcal W$.
We denote by $W>  0$  $(W \geq 0)$, that a matrix or matrix function $W$ is symmetric and positive definite (semi-definite).

The general class of differential algebraic equations (DAEs) (or descriptor systems as they are called in the control context), consists of implicitly defined equations of the form
\begin{eqnarray}
	\label{eqn:descriptorSystem}
	F(t,x(t),\dot{x}(t),u(t)) &=& 0,\\
		y(t) - G(t,x(t),u(t)) &=& 0,\nonumber
\end{eqnarray}
on some time interval $\mathbb T= [t_0,t_{f}]$
with
\begin{align*}
    F\colon \mathbb T\times\mathbb{D}_{x}\times\mathbb{D}_{\dot{x}}\times\mathbb{D}_{u}\to \mathbb R^\ell\qquad\text{and}\qquad
    G\colon\mathbb T\times\mathbb{D}_{x}\times\mathbb{D}_{u}\to \mathbb R^p,
\end{align*}
with open domains, vector spaces, or manifolds
$\mathbb{D}_{x}$, $\mathbb{D}_{\dot{x}}$, $\mathbb{D}_{u}$.
We refer to $x$, $u$, and $y$, as the \emph{state}, \emph{input}, and \emph{output}, respectively.

In industrial practice, the models typically have (control) parameters and/or may have uncertain components such as e.g, unmodeled quantities, uncertainty in parameters, or disturbances, but we will not discuss these topics here and we will also not discuss the typical control tasks, such as feedback control, or optimal control, see \cite{MehU23} for a survey of the control aspects. In the following we therefore assume that the input function is either used for the interconnection of submodels or is given as a fixed function.

{
We consider solutions in the spaces of continuous or continuously differentiable functions. One may, however,  also use weaker solution concepts, see e.g.  \cite{KunM24,RabR96a,RabR96b,Tre13}.
For  automatic modelling in industrial practice, one has to consider over- and under-determined systems, in particular,  to allow redundancies (such as control parameters) in the model.

DAEs have been studied since the early 1960s and they form the basis of most modern simulation tools, see e.g. \cite{BreCP96,KunM24,LamMT13,Ria08} for  monographs covering the topic.
  One  of the major difficulties when modeling with DAEs is that the solution may require stronger  differentiability assumptions on $F$ than in the case of ordinary differential equations, since the algebraic equations may lead to so called \emph{hidden constraints} depending on higher order derivatives that are not directly visible in the equations. The constraints, and in particular the hidden constraints, then also restrict the set of potential initial conditions. To address these difficulties, several different \emph{index concepts} have been developed that allow to characterize the extra differentiability (regularity) requirements and the set of consistent initial conditions, see \cite{Meh15} for a survey and comparison.

For a given input one obtains the DAE
\begin{equation}
    \label{eq:rDAE}
    \mathcal{F}(t,x(t),\dot{x}(t)) = 0.
\end{equation}
If the Jacobian $\tfrac{\partial}{\partial \dot{x}} \mathcal{F}$ is not square or if it is singular, then a solution $x$, provided such a solution exists, typically  depends on derivatives of $\mathcal{F}$. This means that derivatives of (some of) the equations have to be included in the model to guarantee that in a numerical  simulation these extra constraints are satisfied.

To characterize the level of extra differentiability that is needed to obtain classical solutions, to filter out a complete set of algebraic equations that restrict the possible initial conditions,   and to generalize other index concepts to over- and under-determined systems, the \emph{strangeness index} concept was introduced,  see \cite[Chapters 3 and 4] {KunM24},  for a derivation and detailed analysis. It  is like the more commonly used  \emph{differentiation index}, see~\cite{Cam87a}, based on the \emph{derivative array} of level $\mu$, that is defined as
\begin{equation}
	\label{eq:derivativeArray}
	\widetilde{\mathcal F}_\mu\left(t,x,\eta\right) : = \begin{bmatrix}
		\mathcal{F}(t,x,\dot{x})\\
		\ddt\mathcal{F}(t,x,\dot{x})\\
		\vdots\\
		\left(\ddt\right)^{\!\mu} \mathcal{F}(t,x,\dot{x})
	\end{bmatrix}\in\mathbb{R}^{(\mu+1)\ell}\ \text{with}\ \eta: =\begin{bmatrix}
		\dot{x}\\
		\ddot{x}\\
		\vdots\\
		x^{(\mu+1)}
\end{bmatrix}\in \mathbb R^{(\mu+1)n}.
\end{equation}
To give a brief overview of analysis, we  make several assumptions, see \cite[Chapter 4]{KunM24}  and \cite[Section 2]{MehU23}. We first assume that  the solution set  of \eqref{eq:derivativeArray}
%
\begin{equation}
    \label{eqn:nonlinDAE:manifold}
	\mathcal{M}_\mu \vcentcolon= \left\{\left(t,x,\eta\right)\in\mathbb{R}^{(\mu+2)n+1}\ \bigg|\ \widetilde{\mathcal F}_{\mu}\left(t,x,\eta\right) = 0\right\},
\end{equation}
considered as an algebraic equation, is nonempty and (locally) forms a ma\-ni\-fold. Following \cite{KunM98}, one introduces the Jacobians
\begin{align}\label{eqn:nonlinDAE:Jacobians}
	{\mathcal{E}}_\mu &: = \begin{bmatrix}
		\frac{\partial \widetilde{\mathcal F}_{\mu}}{\partial \dot{x}} & \dots & \frac{\partial \widetilde{\mathcal F}_{\mu}}{\partial x^{(\mu+1)}}
	\end{bmatrix}\in\mathbb{R}^{(\mu+1)\ell,(\mu+1)N},\\
{\mathcal{A}}_\mu &: = -\begin{bmatrix}
		\frac{\partial \widetilde{\mathcal F}_{\mu}}{\partial x} & 0 & \dots & 0
	\end{bmatrix}\in\mathbb{R}^{(\mu+1)\ell,(\mu+1)n}.\nonumber
\end{align}
In order to avoid lengthy technical derivations, we assume that on  $\mathcal{M}_\mu$, these Jacobians satisfy
    \begin{equation}
        \label{eqn:nonlinDAE:regularPartAss}
        \rank \begin{bmatrix}
            \mathcal{E}_\mu & \mathcal{A}_\mu
        \end{bmatrix} = r,
    \end{equation}
and
    \begin{equation}
        \corank \begin{bmatrix}
            \mathcal{E}_\mu & \mathcal{A}_\mu
        \end{bmatrix} - \corank \begin{bmatrix}
            \mathcal{E}_{\mu-1} & \mathcal{A}_{\mu-1}
        \end{bmatrix} = v,
    \end{equation}
where  $v$  counts the number of redundancies (equations of the form $0=0$ in the system).
The derivative array then allows to extract the algebraic equations, and thus characterizes the required level of differentiability and defines the set of consistent initial values.  For this, we further assume that the matrix $\mathcal{E}_{\mu}$ defined in~\eqref{eqn:nonlinDAE:Jacobians} satisfies
    \begin{equation}
        \rank \mathcal{E}_{\mu} = r-a\qquad\text{on $\mathcal{M}_\mu$}.
    \end{equation}

Under this assumption (together with a theorem on the existence locally smooth full rank decompositions, \cite[Thm.~4.1.3]{KunM24}), there exists (locally)  a smooth matrix function $Z_{a}\colon \mathcal{M}_{\mu}\to \R^{(\mu+1)\ell,a}$ with pointwise maximal rank on $\mathcal{M}_\mu$ that satisfies
\begin{equation}
    \label{eqn:nonlinDAE:hyp:algebraicSelector}
    Z_{a}^T \mathcal{E}_\mu = 0.
\end{equation}
The (linearized) algebraic equations are then encoded in the matrix function
\begin{equation}
    \label{eqn:nonlinDAE:hyp:linearizedAlgebraicEquations}
    Z_{{a}}^T \tfrac{\partial \widetilde{\mathcal{F}}_{\mu}}{\partial x}.
\end{equation}
To ensure that these algebraic equations are solvable for $a$ unknowns,  one  shows that the matrix in~\eqref{eqn:nonlinDAE:hyp:linearizedAlgebraicEquations} has full rank, since then \eqref{eqn:nonlinDAE:hyp:algebraicSelector}  implies that
\begin{displaymath}
    \rank Z_{{a}}^T \mathcal{A}_{\mu} = \rank Z_{{a}} \tfrac{\partial \widetilde{\mathcal{F}}_{\mu}}{\partial x} = a.
\end{displaymath}
Again, \cite[Thm.~4.1.3]{KunM24} implies (locally) the existence of a smooth matrix function $T_{a}\colon \mathcal{M}_\mu \to \R^{n,n-a}$ with pointwise maximal rank satisfying
$
    Z_{a}^T \tfrac{\partial \widetilde{\mathcal{F}}_\mu}{\partial x}T_{a} = 0
$
on $\mathcal{M}_\mu$.

The remaining differential equations must be contained in the original DAE \eqref{eq:rDAE} (in contrast to the algebraic equations, which are contained in the derivative array). We then assume further that   $        \rank \tfrac{\partial \mathcal{F}}{\partial \dot{x}}T_{a} = d=\ell-a-v$ on $\mathcal{M}_\mu$
and  again, we employ ~\cite[Thm.~4.1.3]{KunM24} to (locally) obtain a smooth matrix function $Z_{{d}}$ of size $N\times d$ with pointwise maximal rank that satisfies $Z_{{d}}^T \tfrac{\partial \mathcal{F}}{\partial \dot{x}}T_{a} = d$. The matrix function $Z_{d}$ then filters out the ordinary differential equations.

Note that, due to the local character,  all the discussed assumptions only hold in a suitable neighborhood, and it may be necessary to restrict the analysis to subintervals in time and then in a numerical integration switch between these time intervals,  see \cite{KunM18} for a detailed discussion and further analysis.

We summarize the discussed assumptions  in the following assumption see \cite[Hypothesis 4.3.1]{KunM24}.
\begin{assumption}
  \label{hyp:nonLin:nonRegular}
  Consider the extended DAE system \eqref{eq:rDAE}  and associated initial values $(t_0,x^0,\eta^0)$.
There exist integers $\mu$, $r$, $a$, and $v$ such that $\mathcal{M}_\mu$ defined in \eqref{eqn:nonlinDAE:manifold} is nonempty and such that for every $(t_0,x^0,\eta^0)\in\mathcal{M}_\mu$ there exists a (sufficiently small) neighborhood $\mathcal{U}$ in which the following properties hold.
    \begin{enumerate}
       \item[(i)] The set $\mathcal{M}_\mu$ forms a manifold of dimension $(\mu+2)n+1-r$.
        \item[(ii)] We have $\rank\begin{bmatrix}
           \mathcal{E}_\mu & \mathcal{A}_\mu
           \end{bmatrix} = r$ on $\mathcal{M}_\mu\cap \mathcal{U}$.
        \item[(iii)] We have $\corank \begin{bmatrix}
            \mathcal{E}_\mu & \mathcal{A}_\mu
            \end{bmatrix} - \corank \begin{bmatrix}
            \mathcal{E}_{\mu-1} & \mathcal{A}_{\mu-1}
           \end{bmatrix} = v$ on $\mathcal{M}_\mu\cap \mathcal{U}$ (with the convention $\corank\begin{bmatrix}
              \mathcal{E}_{-1} & \mathcal{A}_{-1}
           \end{bmatrix} = 0$).
        \item[(iv)] We have $\rank \mathcal{E}_\mu = r-a$ on $\mathcal{M}_\mu\cap\mathcal{U}$, such that there exist locally smooth matrix functions $Z_{a}$ and $T_{a}$ of size $(\mu+1)\ell\times a$ and $n\times (n-a)$, respectively, and pointwise maximal rank, satisfying $Z_{a}^T\mathcal{E}_\mu = 0$, $\rank Z_{{a}}^T\mathcal{A}_\mu = a$, and $Z_{a}^T\tfrac{\partial \widetilde{F}_\mu}{\partial x}T_{a} = 0$ on $\mathcal{M}_\mu\cap\mathcal{U}$.
        \item[(v)] We have $\rank \tfrac{\partial \mathcal{F}}{\partial \dot{x}}T_{{a}} = d : = \ell-a-v$ on $\mathcal{M}_\mu\cap\mathcal{U}$ such that there exists a smooth matrix function $Z_{d}$ of size $n\times d$ and pointwise maximal rank, satisfying $\rank Z_{d}^T \tfrac{\partial \widetilde{F}}{\partial \dot{x}}T_{a} = d$.
    \end{enumerate}
\end{assumption}
The smallest value $\mu$ such that $\mathcal{F}$ satisfies Assumption~\ref{hyp:nonLin:nonRegular}, is called the \emph{strangeness index} of the DAE and it generalizes the more commonly used \emph{differentiation index} $\nu$ to under- and over-determined system. If $\mu=0$, then the DAE is called \emph{strangeness-free}, and one can use the projectors $Z_a$ and $Z_d$ to (locally) construct from the derivative array a \emph{reduced DAE}
\begin{equation}
    \label{eqn:nonlinDAE:sfree}
    \widehat{\mathcal{F}}(t,x,\dot{x}) : = \begin{bmatrix}
    \widehat{\mathcal{F}}_{d}(t,x,\dot{x})\\
    \widehat{\mathcal{F}}_{a}(t,x)
    \end{bmatrix}= \begin{bmatrix}
   (Z_{d}^T \mathcal{F})(t,x,\dot{x})\\
    (Z_{a}^T \widetilde{\mathcal{F}}_\mu)(t,x)
    \end{bmatrix}=0.
\end{equation}
It is shown in~\cite{KunM01} that the reduced quantities $\widehat{\mathcal{F}}_{a}$ and $\widehat{\mathcal{F}}_{d}$ are independent of higher derivatives of~$x$.

For a given input $u$, the reduced DAE \eqref{eqn:nonlinDAE:sfree} is  ideally suited for numerical simulation, since the second equation explicitly displays all the algebraic constraints and thus the set of consistent initial conditions. Moreover, for the time-discretization of the first equation one can use any implicit numerical integration method that is also good for ordinary differential equations.

While the strangeness index is the most general index concept, it reduces to the differentiation index, which is equal to the strangeness index plus one, if the system is not over- or under-determined and not an ordinary differential equation. In many applications (circuit simulation, multibody dynamics, (in)compressible flow, power systems) the reduced  DAE \eqref{eqn:nonlinDAE:sfree} is trivially obtained from the structure of the equations, see \cite[Section 4.2]{KunM24} for several applications and for general DAEs this is also the case if all the algebraic constraints are explicitly present in the model. A particularly nice set are semi-explicit DAEs,
$\dot x_1=f(x_1,x_2)$, $0=g(x_1,x_2)$, where  the derivative array can be simplified by only differentiating the second equation, and where the system is strangeness-free whenever the
Jacobian $\frac{\partial g}{\partial x_2}$ is invertible, since then the second equation can be solved at any time point.

For the analytical treatment of the initial value problem, for  a given control function $u$, and fixed initial data $(t_0,x^0,\eta^0)\in{\mathcal M}_\mu$,  one  can
locally parameterize $\mathcal M_\mu$  by $(\mu+2)n+1-r$ parameters chosen from $(t,x,\eta)$ in such a way
that discarding the associated columns from
$\mathcal F_{\mu;t,x,\eta}(t_0,x^0,\eta^0)$
does not lead to a rank drop.

Without restriction, one then writes~$x$ as $(x_1,x_2,x_3)$
with $x_1\in{\mathbb R}^d$, $x_2\in{\mathbb R}^{n-a-d}$, $x_3\in{\mathbb R}^a$,
one chooses $(x_1,x_2)$ as further parameters and assumes that $\eta=(\eta_1,\eta_2)$
with $\eta_1$ as the chosen parameters. Then, there (locally) exist functions ${\mathcal G}_3$, corresponding to~$x_3$,
and ${\mathcal G}_\eta$, corresponding to $\eta$ such that
\[
\mathcal F_\mu(t,x_1,x_2,{\mathcal G}_3(t,x_1,x_2,\eta_1),{\mathcal G}_\eta(t,x_1,x_2,\eta_1))\equiv0.
\]
Moreover,  there exists a function~${\mathcal R}$ such that
\[
x_3={\mathcal G}_3(t,x_1,x_2,\eta_1)={\mathcal R}(t,x_1,x_2)
\]
giving
\begin{equation}\label{ift1ouq}
F_\mu(t,x_1,x_2,{\mathcal R}(t,x_1,x_2),{\mathcal G}_\eta(t,x_1,x_2,\eta_1))\equiv0.
\end{equation}
Similarly, we can choose~$T_2$ of
Assumption~\ref{hyp:nonLin:nonRegular} as
\[
T_2(t,x_1,x_2)=
 \begin{bmatrix} I\\{\cal R}_{x_1,x_2}(t,x_1,x_2)\end{bmatrix},
\]
and  there a matrix function~$Z_1$
which only depends on the original variables $(t,x,\dot x)$ giving the reduced DAE
\begin{equation}\label{eq:redqou}
\hat {\mathcal F}(t,x,\dot x)=
 \begin{bmatrix} \hat {\mathcal F}_1(t,x,\dot x)\\\hat {\mathcal F}_2(t,x)\end{bmatrix} =0,
\end{equation}
with
\[
\begin{split}
\hat F_1(t,x_1,x_2,x_3,\dot x_1,\dot x_2,\dot x_3)&=
Z_1^TF(t,x_1,x_2,x_3,\dot x_1,\dot x_2,\dot x_3),
\\
\hat F_2(t,x_1,x_2,x_3)&=Z_2^TF_\mu(t,x_1,x_2,x_3,{\mathcal G}_\eta(t,x_1,x_2,\eta_{1}^0)).
\end{split}
\]
We then have the following existence and uniqueness result, see \cite[Theorems 4.1.14 and  4.1.15]{KunM24}.
\begin{theorem}
\label{th:nlsuffou}
Let $\mathcal F$ as in \eqref{eq:rDAE} be sufficiently smooth and satisfy
Assumption~\ref{hyp:nonLin:nonRegular}
with characteristic values $\mu$, $a$, $d$, $v$
as well as with characteristic values $\mu+1$, replacing~$\mu$, and the same values for~$a$, $d$, $v$.
Let $(t_0,x^0,y^0)\in{\mathcal M}_{\mu+1}$ be given
and let the parameterization~$\eta_1$ in $F_{\mu+1}$ include~$\dot x_2$.
Then, for every function $x_2\in C^1({\mathbb I},{\mathbb R}^{n-a-d})$
with $x_2(t_0)=x_{2}^0$, $\dot x_2(t_0)=\dot x_{2}^0$, the reduced DAE
(\ref{eq:redqou}) has a  unique
solution~$x_1$ and~$x_3$ satisfying $x_1(t_0)=x_{1}^0$.
Moreover, the so obtained function $x=(x_1,x_2,x_3)$
locally solves the original problem.

Furthermore, every sufficiently sufficiently differentiable solution of the original DAE also solves the strangeness-free DAE~\eqref{eqn:nonlinDAE:sfree}.
\end{theorem}

Note that this result is stated for a given, sufficiently smooth, control function $u$. In the context of  control problems, the analysis is more intricate, see \cite[Sections 3.6, 3.7 4.5, 4.6]{KunM24} or the survey \cite{MehU23}.

The  summary of the theory discussed so far provides the basis for the analysis of general DAE systems,  and how to employ the reduced DAE in numerical simulations and control, see \cite{KunM24} for the full theory. Most numerical integration methods for the initial value problem are invariant under the exact choice of the projectors $Z_d, Z_a$ and thus it is not necessary to choose them in a smooth way.

However, most of the challenges that arise in the coupling of DAE systems from different physical domains remain.
So  one may ask whether there exists a general model class where the assumptions are easily verified, and where the numerical methods are such that  automated modeling becomes feasible across different physical domains and scales. Such a class is discussed in the next section.}

\section{Port-Hamiltonian differential-algebraic systems}\label{sec:phdae}
{In order to overcome some of the difficulties with DAEs discussed in the previous section, in particular the multiphysics modelling,} we now  introduce  energy based modeling via the model class of nonlinear (dissipative) \emph{port-Hamiltonian DAE (pHDAE) systems} and its  many important  properties, see
\cite{MehM19,MehU23,Mor24} for  details. The main motivation for using this class is that the only universal quantity for the coupling of systems from different physical domains is the transfer of energy between the different components. In the last $25$  years this motivation has lead to a large increase in the research on energy based modeling, see \cite{RasCSS20,SchJ14}  for  surveys. We employ the following definition from \cite{MehM19}.
\begin{definition}
	\label{def:pHDAE}
	Consider a time interval $\mathbb T$, a state space $\mathcal X$ of functions $x: \mathbb T \to \mathbb R^n$, and an extended space $\mathcal S := \mathbb T\times\mathcal X$.  Then a \emph{port-Hamiltonian descriptor system} (pHDAE) is of the form
	\begin{eqnarray}
			E(t,x)\dot x + \rho (t,x) &=& (J(t,x)-R(t,x))\eta(t, x) + (B(t,x) - P(t,x))u, \nonumber\\
			y &=& (G(t,x) + P(t,x))^T \eta(t, x) + (S(t,x) - N(t,x))u,	\nonumber\\	\label{eqn:pHDAE}
		\end{eqnarray}
	with state $x\colon \mathcal S \to\mathbb R^n$,  input $u\colon \mathcal S \to\mathbb R^m$,  and output $y\colon \mathcal S \to\mathbb R^m$, where
	\begin{align*}
		\rho,\eta&\in C(\mathcal S ,\mathbb R^{\ell}),
		& E&\in C(\mathcal S ,\mathbb R^{\ell, n}),
		& J,R&\in C(\mathcal S ,\mathbb R^{\ell,\ell})\\
		G,P&\in C(\mathcal S ,\mathbb R^{\ell, m}),
		& S,N&\in C(\mathcal S ,\mathbb R^{m, m}),
	\end{align*}
together with an associated function $\mathcal H \in C^1(\mathcal S ,\mathbb R)$, called the \emph{Hamiltonian} of~\eqref{eqn:pHDAE}.  Furthermore, the following properties are required to hold:
	 \begin{enumerate}
	 	\item[(i)] The matrix functions
	 		\begin{eqnarray}
	 		\label{eqn:pHDAE:prop1}
	 		\Gamma &:=& \begin{bmatrix}
	 				J & B\\
	 				-B^T & N
	 			\end{bmatrix}\in C(\mathcal S ,\mathbb R^{(\ell+m),(\ell+m)}),\\
	 			W&:=& \begin{bmatrix}
	 				R & P\\
	 				P^T & S
	 			\end{bmatrix}\in C(\mathcal S ,\mathbb R^{(\ell+m),(\ell+m)}),\nonumber
	 		\end{eqnarray}
	 		satisfy $\Gamma= -\Gamma^T$ and $W = W^T \geq 0$ in $\mathcal S$.
	 	\item[(ii)] The Hamiltonian satisfies
	 		\begin{equation}
	 			\label{eqn:energyCondition}
\tfrac{\partial}{\partial {x}}\mathcal H (t,x) = E^T(t,x)\eta(t,x) \quad\text{and}\quad \tfrac{\partial}{\partial t} \mathcal H (t,x) = \eta^T(t,x) \rho(t,x)
	 		\end{equation}
	 		in $\mathcal S$ along any solution of~\eqref{eqn:pHDAE}.
	 \end{enumerate}
	 If the pHDAE has no inputs and the output equation is omitted, then \eqref{eqn:pHDAE} is called a \emph{(dissipative) Hamiltonian DAE (dHDAE)}.
\end{definition}		
Definition~\ref{def:pHDAE} reduces to the classical representation of ordinary  pH systems if $E $ is the  identity matrix,  and  $\rho\equiv 0$, see e.g. \cite{SchJ14}, where pH systems are defined via a Dirac structure.
In many applications the Hamiltonian is  \emph{nonnegative}, i.e.
$\mathcal H (t,x(t)) \geq 0$ for all $(t,x)\in \mathcal S$ that are solutions of~\eqref{eqn:pHDAE} but for the following results it is sufficient that it is bounded from below. 
%
\begin{theorem}\label{thm:pbe}
A pHDAE of the form \eqref{eqn:pHDAE} satisfies the \emph{power balance equation (PBE)}
  \begin{equation}\label{eq:powerBalanceEq}
    \dd{}{t}\mathcal H(t,x) = - \bmat{x \\ u}^TW\bmat{x \\ u} + y^Tu=:\Delta (x,t)+ \Sigma(x,t)
  \end{equation}
  along any solution $x$ and for any input $u$.
  In particular, the \emph{dissipation inequality}
  \begin{equation}\label{eq:dissIneq}
    \mc H(t_2,x(t_2)) - \mc H(t_1,x(t_1)) \leq \int_{t_1}^{t_2}y(\tau)^Tu(\tau)\td\tau
  \end{equation}
  holds.
\end{theorem}
%
The PBE states that the internally stored energy, represented by the Hamiltonian, only changes by dissipation, represented via the dissipation function $\Delta(x,t)$, and the supplied energy, represented by the function $\Sigma(x,t)$.
Theorem~\ref{thm:pbe} also implies that the Hamiltonian $\mc H$ is a candidate for a Lyapunov function and if it is, then the system is \emph{Lyapunov stable}.  Furthermore,  pHDAE systems are \emph{passive}, see e.g. \cite{MehU23}. This, in particular, means that internally they do not generate energy.

One also has invariance of the class under  state space transformations, see \cite{MehM19}.
{
\begin{theorem}\label{thm:varTrans}
  Consider a pHDAE of the form \eqref{eqn:pHDAE}. Let $\tilde{\mc X}$ be a second state space, let $\tilde{\mc S}:=\mB I\times\tilde{\mc X}$, let $x=\varphi(t,\tilde x)\in\mc C^1(\tilde{\mc S},\mc X)$ be a local diffeomorphism (with respect to $\tilde x$) and let $U\in\mc C(\tilde{\mc S},\R^{\ell,\ell})$ be pointwise invertible.
  Then for $\tilde E=U^T(E\circ\varphi)\pa_{\tilde x}\varphi$, $\tilde J=U^T(J\circ\varphi)U$, $\tilde R=U^T(R\circ\varphi)U$, $\tilde B=U^T(B\circ\varphi)$, $\tilde P=U^T(P\circ\varphi)$, $\tilde \eta=U^{-1}(\eta \circ\varphi)$ and $\tilde \rho=U^T(\rho\circ\varphi+(E\circ\varphi)\pa_t\varphi)$,
  where we set $(F\circ\varphi)(t,\tilde x)=F(t,\varphi(t,\tilde x))$ for any $F\in\mc C(\mc S,\cdot)$, and  $\tilde{\mc H}(t,\tilde x):=(\mc H\circ\varphi)(t,\tilde x)$, the system
  \begin{eqnarray}\nonumber
    \tilde E\dot{\tilde x} + \tilde \rho &=& (\tilde J-\tilde R)\tilde \eta + (\tilde B-\tilde P)u, \\
    y &=& (\tilde B+\tilde P)^T\tilde \eta + (S-N)u,\label{eq:tformPHDAE}
  \end{eqnarray}
  is a pHDAE with Hamiltonian function $\tilde{\mc H}$.
  Furthermore,  to any solution $(\tilde x,u,y)$ of \eqref{eq:tformPHDAE} there corresponds a solution $(x,u,y)$ of \eqref{eqn:pHDAE} with $x(t)=\varphi(t,\tilde x(t))$ and  if $\varphi(t,\cdot)$ is a global diffeomorphism for all $t\in\mB I$, then the two systems are equivalent.
\end{theorem}
 Theorem~\ref{thm:varTrans} implies, in particular, that we can choose appropriately rescaled state variables and representations without changing the input-output map and the power balance equations and just representing the Hamiltonian in different variables. This turns out to be a key factor in simulation of coupled systems across different physical domains.}

{Another property that is relevant for numerical time integration is that one can easily  make a  system with the properties \eqref{eqn:pHDAE}
\emph{autonomous} by introducing
$\widehat{x} := [x^T, t]^T$ and  obtaining an extended system of the form~ \eqref{eqn:pHDAE}. Note that in this extended system the term $\rho(t,x)$ vanishes.

An important property that makes  pHDAEs ideal for automated network based modular modeling is that the class is invariant under power conserving interconnection. Let us consider two autonomous pHDAEs (leaving out the arguments for ease of presentation)
\begin{align*}
  E_i\dot x_i &= (J_i-R_i)\eta_i + (B_i-P_i)u_i, \\
  y_i &= (B_i+P_i)^T\eta _i + (S_i-N_i)u_i,
\end{align*}
with Hamiltonian $\mc H_i$, for $i=1,2$, and assume that the aggregated input $u=(u_1,u_2)$ and output $y=(y_1,y_2)$ satisfy a linear interconnection relation $Mu+Ny=0$ for some $M,N\in\fc*{\R^{k,m}}$. Then the aggregated system can be written as a pHDAE of the form
\begin{align*}
  \bmat{E & 0 & 0 \\ 0 & 0 & 0 \\ 0 & 0 & 0 \\ 0 & 0 & 0}\bmat{\dot x \\ \dot{\hat u} \\ \dot{\hat y}} &=
  \bmat{\Gamma-W & \mat{0 & 0 \\ I_m & -M^T}\hskip-5pt \\ \mat{0 & -I_m \\ 0 & M}\hskip-5pt & \mat{0 & -N^T \\ N & 0}\hskip-5pt}
  \bmat{\eta \\ \hat u \\ \hat y \\ 0}
  + \bmat{0 \\ 0 \\ I_m \\ 0}u, \\
  y &= \hat y,
\end{align*}%
with Hamiltonian $\mc H=\mc H_1+\mc H_2$, with the new variables $\hat u,\hat y\in\R^n$, $x=(x_1,x_2)$, $\eta=(\eta_1,\eta_2)$, $m=m_1+m_2$,   $E=\diag(E_1,E_2)$, $\Gamma=\Pi\diag(\Gamma_1,\Gamma_2)\Pi^T$ and $W=\Pi\diag(W_1,W_2)\Pi^T$, and where $\Pi\in\R^{\ell+m,\ell+m}$ is a permutation matrix, see \cite{MehM19,Mor24} for a detailed analysis.

Due to the input-output structure, the class of pHDAEs allows that state space models of different type and dimension can be combined. These could be infinite or finite dimensional   systems obtained by the space discretization of partial differential equations, or surrogate models obtained by model order reduction or data based modeling, \cite{AntBG20,BenCOW17}.  This makes the class of pHDAEs  particularly useful in hardware-in-the loop approaches, see e.g. \cite{Poe22}, and space-time-model adaptive simulation and optimization methods, see Section~\ref{sec:adaptive}. This also means that one can analyze each subsystem separately concerning its properties, in particular with respect to regularity and the index of pHDAEs.

In many circumstances, e.g. in linear pHDAE systems  or systems with quadratic and convex Hamiltonian,  it has been shown that the strange\-ness index of a pHDAE is at most one and that the system is regular if and only the coefficients $E,J,R$ do not have a common nullspace, \cite{MehMW18,MehS23}.

Furthermore, variables that would have to be differentiated to obtain the strangeness-free reduced model  do not contribute to the Hamiltonian, see \cite{MehS23} for a detailed analysis via normal forms. Under some further conditions on the Hamiltonian it can be shown that similar results also hold in the general case, see \cite{Mor24}.
 This implies, in particular, that after interconnection of two pHDAE systems of the same strangeness index the index cannot increase.

As a consequence many of the numerical and analytical difficulties associated with strangeness index of two or higher do not arise in automated network based pHDAE modeling.

Finally, if one uses appropriate structured (Petrov-)Galerkin projections to restrict the state space, as in space discretization via finite element, finite difference, finite volume methods, or in model reduction, then one stays in the class of pHDAE systems, see e.g. \cite{PolS10,Mor24,SerMH19}. This follows directly from the fact that the symmetry or skew-symmetry struture of the coefficient functions is preserved by structured projections,
see \cite{MehS23} for a detailed characterization of structure preserving transformations that reflect the geometric and algebraic structure and associated normal forms in the linear case.

To illustrate the generality and wide applicability of pHDAE systems in multi-physics modeling  consider the following two example classes.
\begin{example}\label{ex:circuit}{\rm
The first domain, where automated network based DAE modeling was applied with great success was the field of  circuit simulation. An RLC circuit can be modeled as a directed graph with incidence matrix
$    \calA = \begin{bmatrix} 	\calA_r & \calA_c & \calA_\ell & \calA_v & \calA_i \end{bmatrix}$,
partitioned into components associated with resistors, capacitors,  inductors, voltage sources, and current sources; see e.g.~\cite{Fre11} for further details. Let $V$ denote the vector of voltages at the nodes (except for the ground node at which the voltage is zero). Furthermore, let $I_\ell$, $I_v$, and $I_i$ denote the vectors of currents along the edges for the inductors, voltage sources, and current sources, respectively, while $V_v$ and $V_i$ denote the vectors of voltages across the edges for the voltage sources and current sources. Considering, for simplicity of presentation, only linear devices and using Kirchhoff's current and voltage law combined with the so-called branch constitutive relations  yields a pHDAE of the form \eqref{eqn:pHDAE} given by
\begin{displaymath}
	\resizebox{.95\linewidth}{!}{$\begin{aligned}
	\begin{bmatrix}
		\calA_c\mathsf{C}\calA_c^T & 0 & 0\\
		0 & \mathsf{L} & 0\\
		0 & 0 & 0
	\end{bmatrix}\begin{bmatrix}
		\dot{V}\\
		\dot{I}_{\ell}\\
		\dot{I}_v
	\end{bmatrix} &= \begin{bmatrix}
		-\calA_r \mathsf{R}^{-1} \calA_r^T & -\calA_\ell & -\calA_v\\
		\calA_\ell & 0 & 0\\
		\calA_v & 0 & 0
	\end{bmatrix}\begin{bmatrix}
		V\\
		I_{\ell}\\
		I_v
	\end{bmatrix} + \begin{bmatrix}
		-\calA_i & 0\\
		0 & 0\\
		0 & -I
	\end{bmatrix}\begin{bmatrix}
		-I_i\\
		V_v
	\end{bmatrix},\\
	\begin{bmatrix}
		V_i\\
		-I_v
	\end{bmatrix} &=  \begin{bmatrix}
		\calA_i & 0\\
		0 & 0\\
		0 & -I
	\end{bmatrix}^T \begin{bmatrix}
		V\\
		I_{\ell}\\
		I_v
	\end{bmatrix}
	\end{aligned}$}
\end{displaymath}
and Hamiltonian
\begin{displaymath}
	\hamiltonian(V,I_\ell) = V^T \calA_c\mathsf{C}\calA_c^T V + I_\ell^T \mathsf{L} I_\ell,
\end{displaymath}
 associated with the stored energy in the capacitors and inductors.  Note that the current $I_r$ is not present in the Hamiltonian. The symmetric positive definite matrices $\mathsf{R}$, $\mathsf{C}$, and $\mathsf{L}$ are defined via the defining properties of the resistors, capacitors, and inductors.
Note that in this case, the pHDAE is linear with a quadratic Hamiltonian, and no feedthrough term.
 For nonlinear device elements a corresponding relation is easily obtained, see \cite{Mor24}.
For more general circuits, recently in \cite{NedPS22,ShaCE22} new energy based formulations have been suggested. As discussed before, the success of pHDAE modeling is particularly visible in the coupling of circuit and magnetic field models \cite{AltGPSS26}. Another recent development is the formulation of dynamic iteration schemes for coupled pHDAE systems and their use in circuit simulation in \cite{GunBJR21}.
}
\end{example}

\begin{example}\label{ex:mbs}{\rm
    A second classical domain of applications arises in constrained multibody systems, see e.g. \cite{EicF98,KunM24,RabR00,Sim13}. These can be easily expressed as pHDAE system \cite{Sch13,SchJ14}. We demonstrate this for linear multibody systems  which, after adding a tracking output, lead to  a system
\begin{equation}
    \label{eq:mblin}
    \begin{aligned}
	M \ddot{p}  &= - D\dot{p} - Kp -G^T \lambda+ B_1 u, \\
	0 & = Gp ,\\
y& = B_1 \dot p,
    \end{aligned}
\end{equation}
where~$p$ is a vector of positions or position coordinates in a finite element model, $v=\dot p$ is a vector of the associated  velocities, $M$~is the mass matrix, which is usually  symmetric and positive definite,
$Gp=0$~describes the constraints and~$\lambda$ is the associated vector of
Lagrange multipliers that penalize the violation of the constraints. Analogous structures arise in the nonlinear case.

If no redundant constraints occur in \eqref{eq:mblin},
then the system has strangeness index~$\mu=2$ (differentiation index~$\nu=3$), see \cite{EicF98,KunM24}.
For this reason it is typically necessary to use a regularization procedure to make the system well suited for numerical simulation and control, see e.g. \cite{EicF98,KunM24,Sim13}. One possibility is to replace the original constraint by its time derivative $0=G \dot p=Gv$. By transforming to first order form with state
$   x = \begin{bmatrix} v^\top & p^\top & \lambda^\top  \end{bmatrix}^\top$,
one directly obtains a linear  pHDAE of the form \eqref{eqn:pHDAE} with $P =0$, $S-N=0$, and
\[
    E = \begin{bmatrix} M & 0 & 0 \\ 0 & K & 0 \\ 0 & 0 & 0 \end{bmatrix}, \
    R = \begin{bmatrix} D &  0 & 0 \\
0 & 0 & 0 \\ 0 & 0  & 0 \end{bmatrix}, \
    B = \begin{bmatrix} B_1 \\ 0 \\ 0 \end{bmatrix}, \
    J = \begin{bmatrix} 0 & -K &  G^T \\ K & 0 & 0 \\ -G & 0 & 0 \end{bmatrix},
\]
and  quadratic Hamiltonian $    \mathcal H(x)= \frac{1}{2} (v^\top M v +p^\top  K  p)$.
Note that the Lagrange multiplier does not contribute to the Hamiltonian.
}
\end{example}
It has become a well-established modeling paradigm that many physical domains either directly or after some simple modifications lead to pHDAE systems. {For many further applications, see  \cite{JacZ12,MehU23,RasCSS20,SchJ14}. We do not discuss the control aspects here, but recent research shows again the advantages of the pHDAE structure. As an example one may look at the hard and partly open problem of output-feedback stabilization for DAEs, which turns out to be quite simple for pHDAEs, \cite{ChuM25b,ChuM25a}.

In the next section we discuss two applications in gas networks and district heating networks where another advantage of the systems theoretic network based modeling approach with pHDAEs is demonstrated, the fact that one can employ model catalogs for model adaptivity.
\section{Model catalogs and their use in simulation and optimization}\label{sec:adaptive}
In this section we discuss (pHDAE) model catalogs and their use in the system theoretic setting.} Consider the state representation  in the input-output map , see  Figure~\ref{fig:bbmodel}, with a black-box which represents a full catalog of models as e.g. displayed in Figure~\ref{fig:Hiearchy}.
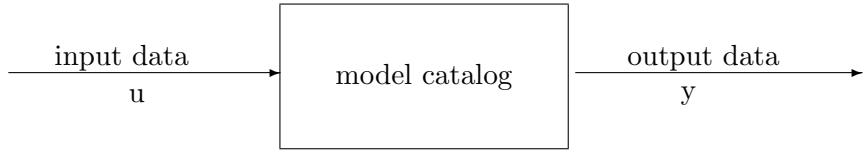
\begin{figure}[H]
\unitlength=1.00mm
\special{em:linewidth 0.4pt}
\linethickness{0.4pt}
\vskip .5truecm
\hskip -2 truecm
\begin{picture}(131.00,24.00)
\put(18.00,15.00){\vector(1,0){36.00}}
\put(93.00,15.00){\vector(1,0){38.00}}
\put(33.00,17.00){\makebox(0,0)[cc]{input data}}
\put(110.00,17.00){\makebox(0,0)[cc]{output data}}
\put(35.00,12.00){\makebox(0,0)[cc]{u}}
\put(108.00,12.00){\makebox(0,0)[cc]{y}}
\put(54.00,5.00){\framebox(38.00,19.00)[cc]{model catalog}}
\end{picture}
\caption{pHDAE black box model}
        \label{fig:bbmodel}
\end{figure}

\begin{figure}[H]
        \centering
        \includegraphics[scale=0.3]{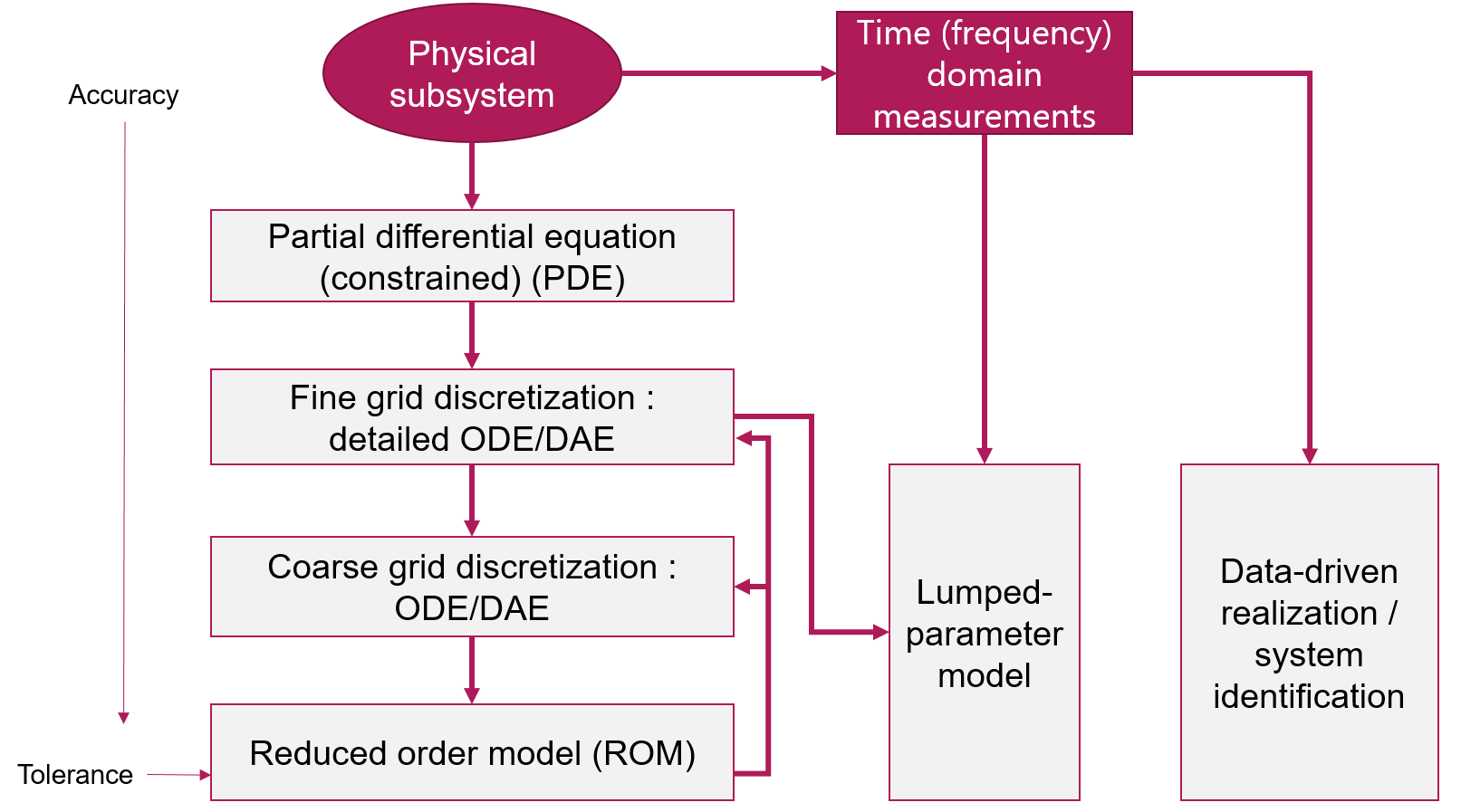}
        \caption{Model catalog.}
        \label{fig:Hiearchy}
\end{figure}
{
Such model catalogs were derived in the context of digital twin modeling for electrical machines in \cite{CheMS22,CheSMGL26}, in \cite{DomHLMMT21} for gas transport networks, see also the upcoming book \cite{MarLHLPTU26}, in \cite{DaeMRS24,HauMMMMRS20,KruMS21} for the simulation and optimization of district heating networks, or in \cite{Poe22} for multiphysics simulation in the automotive sector.

These catalogs are particularly successful in real time simulation and optimization of large scale systems, because they allow adaptivity in model, space discretization and time discretization. In this way it is possible to generate a trade-off between accuracy and computational efficiency   by switching the model as well as the space-time  discretization mesh using an adaptive strategy that is based on error and sensitivity estimates. This allows to achieve a guaranteed error tolerance  for  problems, where the finest and most accurate model is a very good approximation of the physical reality.
}

In the following, we briefly illustrate {this approach for  gas transport networks, see e.g. \cite{DomDSLM18,MehSS18,StoM18} as well as  district heating networks see \cite{DaeMRS24}.  For the flow in the pipes detailed port-Hamiltonian model catalogs based on the compressible, respectively incompressible, 1D Euler equations of fluid dynamics with varying levels of accuracy for the pipe flow, \cite{DomHLMMT21,HauMMMMRS19}.
 These models, which are typically constrained port-Hamiltonian partial differential equations, see \cite{JacZ12,ZwaM24}, lead to pHDAEs after space discretization, respectively model reduction, using projection methods. The other components, such as input sources, valves, compressors and consumers were added in a network based fashion, see Figure~\ref{fig:TestNet3} from \cite{DomDSLM18} for such a simple network.
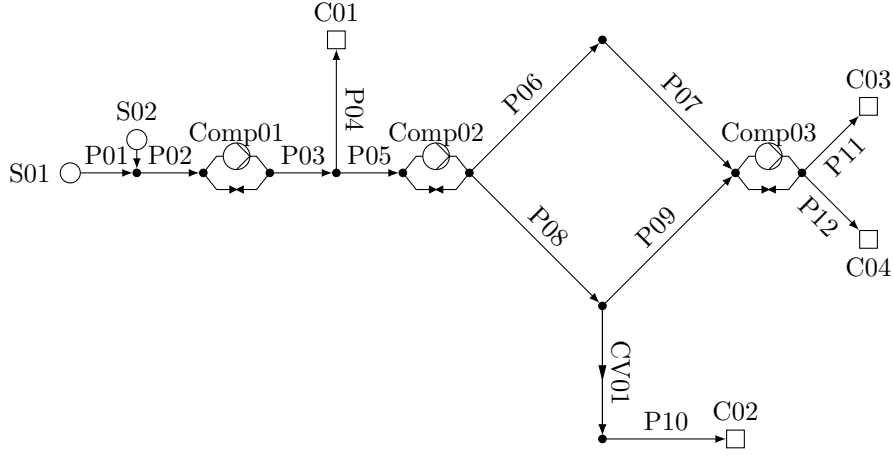
\begin{figure}[!htb]
\centering \small
 \begin{tikzpicture}[scale=0.088]
		\node[circle,fill=black,draw,scale=0.3] (N01) at (20,40) {};
		\node[circle,fill=black,draw,scale=0.3] (N04) at (30,40) {};
		\node[circle,fill=black,draw,scale=0.3] (N05) at (40,40) {};
		\node[circle,fill=black,draw,scale=0.3] (N06) at (50,40) {};
		\node[circle,fill=black,draw,scale=0.3] (N08) at (60,40) {};
		\node[circle,fill=black,draw,scale=0.3] (N09) at (70,40) {};
		\node[circle,fill=black,draw,scale=0.3] (N10) at (90,60) {};
		\node[circle,fill=black,draw,scale=0.3] (N11) at (90,20) {};
		\node[circle,fill=black,draw,scale=0.3] (N12) at (90,0) {};
		\node[circle,fill=black,draw,scale=0.3] (N13) at (110,40) {};
		\node[circle,fill=black,draw,scale=0.3] (N16) at (120,40) {};
    	   	\node[circle,draw,scale=0.8,label={[scale=1]left:{S01}}] (S01) at (10,40) {};
    	   	\node[circle,draw,scale=0.8,label={[scale=1]above:{S02}}] (S02) at (20,45) {};
 		\node[draw,scale=1,label={[scale=1]above:{C01}}] (C01) at (50,60) {};
 		\node[draw,scale=1,label={[scale=1]above:{C02}}] (C02) at (110,0) {};
 		\node[draw,scale=1,label={[scale=1]above:{C03}}] (C03) at (130,50) {};
 		\node[draw,scale=1,label={[scale=1]below:{C04}}] (C04) at (130,30) {};
 		\path[-latex]
 		(S01) edge node[above,sloped,scale=1,text=black] {{P01}} (N01)
 		(N01) edge node[above,sloped,scale=1,text=black] {{P02}} (N04)
 		(N05) edge node[above,sloped,scale=1,text=black] {{P03}} (N06)
 		(N06) edge node[above,sloped,scale=1,text=black] {{P04}} (C01)
 		(N06) edge node[above,sloped,scale=1,text=black] {{P05}} (N08)
 		(N09) edge node[above,sloped,scale=1,text=black] {{P06}} (N10)
 		(N10) edge node[above,sloped,scale=1,text=black] {{P07}} (N13)
 		(N09) edge node[above,sloped,scale=1,text=black] {{P08}} (N11)
 		(N11) edge node[above,sloped,scale=1,text=black] {{P09}} (N13)
 		(N12) edge node[above,sloped,scale=1,text=black] {{P10}} (C02)
 		(N16) edge node[below,sloped,scale=1,text=black] {{P11}} (C03)
 		(N16) edge node[below,sloped,scale=1,text=black] {{P12}} (C04)
 		(S02) edge (N01)
		(N11) edge node[above,sloped] {{CV01}} (N12)
		(N11) edge node[left=4] {} (N12)
 		;
		\draw[fill=black] (N11)++(0,-10)++(0,-1) -- ++(0.5,2) -- ++(-1,0) -- cycle;
		\path[draw] (N04) -- ++(2,2.5) -- node[above=0.1] {{Comp01}} ++(6,0) -- ++(2,-2.5) -- ++(-2,-2.5) --
		node[below=4] {} ++(-6,0) -- ++(-2,2.5);
		\draw[fill=white] (N04)+(5,2.5) circle (2);
		\path[draw] (N04)++(5,4.5) -- ++(2,-2) -- ++(-2,-2);
		\draw[fill=black] (N04)++(5,-2.5)++(-1,0.5) -- ++(2,-1) -- ++(0,1) -- ++(-2,-1) -- cycle;
		\path[draw] (N08) -- ++(2,2.5) -- node[above=0.1] {{Comp02}} ++(6,0) -- ++(2,-2.5) -- ++(-2,-2.5) --
		node[below=4] {} ++(-6,0) -- ++(-2,2.5);
		\draw[fill=white] (N08)+(5,2.5) circle (2);
		\path[draw] (N08)++(5,4.5) -- ++(2,-2) -- ++(-2,-2);
		\draw[fill=black] (N08)++(5,-2.5)++(-1,0.5) -- ++(2,-1) -- ++(0,1) -- ++(-2,-1) -- cycle;
		\path[draw] (N13) -- ++(2,2.5) -- node[above=0.1] {{Comp03}} ++(6,0) -- ++(2,-2.5) -- ++(-2,-2.5) --
		node[below=4] {} ++(-6,0) -- ++(-2,2.5);
		\draw[fill=white] (N13)+(5,2.5) circle (2);
		\path[draw] (N13)++(5,4.5) -- ++(2,-2) -- ++(-2,-2);
		\draw[fill=black] (N13)++(5,-2.5)++(-1,0.5) -- ++(2,-1) -- ++(0,1) -- ++(-2,-1) -- cycle;
 \end{tikzpicture}
\caption{Gas supply network with compressor stations and a control valve.}
\label{fig:TestNet3}
\end{figure}

The network consists of
pipes (P01--P12), sources (S01--S02),
consumers (C01--C04),
compressor stations (Comp01--Comp03)
and one control valve (CV01).
It is modeled
as directed graph $\mathcal{G} \! = \! (\mathcal{J},\mathcal{V})$
with edges~$\mathcal{J}$ and vertices~$\mathcal{V}$.
The set of edges~$\mathcal{J}$
consists of pipes $j \in \mathcal{J}_p$,
compressor stations $c \in \mathcal{J}_c$,
and valves $v \in \mathcal{J}_v$.
Each pipe $j \in \mathcal{J}_p$
is defined as an interval $[x_j^a,x_j^b]$
with a direction from $x_j^a$ to $x_j^b$.
}

Based on the full catalog, in~\cite{Domschke2011b},  a three level catalog was constructed (at that time not yet written in port-Hamiltonian form) for the gas flow simulation based on the \emph{isothermal compressible 1D Euler equations}
consisting  of the continuity and the momentum equation
together with the equation of state for real gases. The top level assumed constant speed of sound and horizontal pipes
which lead to the  model
\begin{align}
\label{M1}
\begin{split}
p_t + \frac{\rho_0 c^2}{A} q_x & = 0, \\
q_t + \frac{A}{\rho_0} p_x + \frac{\rho_0 c^2}{A} \left(\frac{q^2}{p}\right)_{\!\!x} &
= -\frac{\lambda \rho_0 c^2 |q| q}{2DAp},
\end{split}
\end{align}
%
where,
$q = A \rho v/\rho_0$ denotes the mass flow rate
under standard conditions,
$p$ denotes the pressure,
$c = \sqrt{p/\rho}$ the speed of sound,
$A$ the cross-sectional area of the pipe,
$\lambda > 0 $ the Darcy friction coefficient,
$D$ the pipe diameter,
$\rho$ the gas density,
$v$ the gas velocity,
and $\rho_0$ the density
under standard conditions.

Since the velocity of the gas is much smaller than the speed of sound,
as a simplified second level model one  neglects the nonlinear term in the spatial derivative
of the momentum equation in~\eqref{M1}, which leads to the model~\cite{DomschkePhD}
\begin{align}
\label{M2}
\begin{split}
p_t + \frac{\rho_0 c^2}{A} q_x & = 0, \\
q_t + \frac{A}{\rho_0} p_x &
= -\frac{\lambda \rho_0 c^2 |q| q}{2DAp}.
\end{split}
\end{align}
The third level  assumes a stationary state and leads to the  model
\begin{align}
\label{M3}
\begin{split}
q &= \text{const.}, \\
p(x) &= \sqrt{p_{\text{in}}^2 - \frac{\lambda \rho_0^2 c^2 |q| q}{D A^2} x},
\end{split}
\end{align}
where,
$p_{\text{in}} = p(0)$ denotes the pressure
at the inbound of the pipe.
All three models are easily expressed in pHDAE form, and thus the complete model can be built in a modularized form combining different pipeline sections and  models for the other components, see \cite{DomHLMMT21,MarLHLPTU26}.
After space discretization the first two levels are strangeness-free (actually ordinary differential equations, while the third level is also strangeness-free but purely algebraic.

For the space-time discretization, in \cite{DomDSLM18}, an implicit box scheme was used in each pipe using  the three models
and varying discretization step sizes in space and time. Three different refinement strategies based on individual tolerances, maximal error refinement, and maximal error-to cost refinement were studied and it was shown that this leads to savings in CPU time  of $70-80\%$ to achieve an error of $10^{-5}$ for synthetic examples and even more than $90\%$ for the more realistic  test network Figure~\ref{fig:TestNet3}.

{The same  model catalog was also used for the optimization of compressor costs,
 where the model switching was used to guarantee feasible solutions at every optimization step} and it was shown for two industrial examples, documented in the benchmark list \cite{SchABHJKKOPS17}, that a speedup between $13$ and $17$ can be achieved.

A similar approach was used in \cite{DaeMRS24}  to minimize the overall costs in a district heating network to satisfy the heat demand of all the consumers, see \cite{MehSS18}. The considered costs are those for waste incineration,  natural gas to heat the water and for increasing the pressure of the water in the depot. The data were taken from \cite{NusT16} and the water flow  was  modeled
via  1D compressible Euler
equations, given for any edge (pipe) the partial differential equation
  \begin{align*}
    0 &= \frac{\partial \rho_\arc}{\partial t} + \vel_\arc \frac{\partial
        \rho_\arc}{\partial x} + \rho_\arc \frac{\partial \vel_\arc}{\partial
        x},\\
    0 &= \frac{\partial (\rho_\arc \vel_\arc)}{\partial t} + \vel_\arc
        \frac{\partial (\rho_\arc \vel_\arc)}{\partial x} + \frac{\partial \press_\arc}{\partial
        x} \nonumber + \frac{\lambda_\arc}{2 \diam_\arc} \rho_\arc \abs{\vel_\arc}
        \vel_\arc + \grav \rho_\arc h_\arc',
  \end{align*}
where the first equation describes the conservation of mass, whereas the second defines the pressure gradient.
Here the used quantities are the diameter of a pipe~$\arc$, denoted by~$\diam_\arc$, the friction coefficient in the pipe denoted by
$\lambda_\arc$, and $h_\arc'$
denotes the slope of the pipe.
The temperature is denoted by $\temp_\arc$ and  $\heattrans$ and $\soiltemp$ are the heat transfer
coefficient and the soil or pipe wall temperature.
Finally, $g$ is  the gravitational acceleration.

 The equations operate on a directed and connected
graph with two parts, the hot water forward-flow and  the cooled water backward-flow part.
The interconnection of the pipes was performed by using the boundary values at the interconnection points. The network is connected  by the depot where  the cooled water
is heated again, and the consumers who use the temperature
difference to satisfy their thermal energy demand. The system is completed by  energy and  entropy equations and an equation for the incompressibility,
  \begin{align*}
    0 &= \frac{\partial \densintene_\arc}{\partial t} + \vel_\arc
        \frac{\partial \densintene_\arc}{\partial x}  
        +\press_\arc
        \frac{\partial \vel_\arc}{\partial x} - \frac{\lambda_\arc}{2
        \diam_\arc} \rho_\arc \abs{\vel_\arc} \vel_\arc^2 + \frac{4
        \heattrans}{\diam_\arc}(\temp_\arc - \soiltemp),  \\
    0 &= \frac{\partial \densintent_\arc}{\partial t} + \vel_\arc
        \frac{\partial \densintent_\arc}{\partial x} +
        \frac{\lambda_\arc \rho_\arc }{2 \diam_\arc \temp_\arc}
        \abs{\vel_\arc} \vel_\arc^2 + \frac{4 \heattrans}{\diam_\arc}\frac{(\temp_\arc -
        \soiltemp)}{\temp_\arc},\\
  0 &= \frac{\partial \rho_\arc}{\partial t} + \vel_\arc
  \frac{\partial \rho_\arc}{\partial x}.
  \end{align*}
In \cite{HauMMMMRS20} it is shown how to write this complete set of  equations in pH form and it shown that the system is semi-explicit of strangeness-index one due to the incompressibility condition. The index reduction to the strangeness-free formulation is easily obtained by  using the structure, and based on this discretization methods are discussed in \cite{HauM21,Hau24}.

Since current optimization methods are not able to treat the instationary model,  a stationary  model was used together with  $\vel_\arc(x) = \vel_\arc$ and $\rho_\arc(x) = \rho_\arc=\rho$  constant for all pipes. Closing the system by appropriate state equations, these simplifications give the model
  \begin{align*}
    0 &=  \frac{\diff\press_\arc}{\diff x} + \frac{\lambda_\arc}{2
        \diam_\arc} \rho \abs{\vel_\arc} \vel_\arc + \grav
        \rho h_\arc', \\
    0 &=  \vel_\arc \frac{\diff\densintene_\arc}{\diff x} - \frac{\lambda_\arc}{2
        \diam_\arc} \rho \abs{\vel_\arc} \vel_\arc^2 + \frac{4
        \heattrans}{\diam_\arc}(\temp_\arc - \soiltemp), \\
    0 &=  \vel_\arc \frac{\diff\densintent_\arc}{\diff x} +
        \frac{\lambda_\arc \rho }{2 \diam_\arc \temp_\arc}
        \abs{\vel_\arc} \vel_\arc^2 + \frac{4 \heattrans}{\diam_\arc}\frac{(\temp_\arc -
        \soiltemp)}{\temp_\arc},
  \end{align*}
  which  is used as the first and most accurate modeling level in a model catalog. For the second  level, one  neglects the
(small) term $\lambda_\arc / (2 \diam_\arc) \rho \vel_\arc^2
\abs{\vel_\arc}$ and for the third leven one further neglects the term $4
\heattrans / \diam_\arc (\temp_\arc - \soiltemp)$. For the space discretization, the implicit mid-point rule in space was used, and altogether this leads to a highly nonlinear, large scale
finite-dimensional  mathematical program with complementarity constraints.

An adaptive optimization algorithm was then constructed to iteratively solve the
nonlinear program  until one has a feasible solution satisfying a prescribed tolerance.
The algorithm iteratively switches the model level and the grid sizes in the discretization grid for each pipe individually according to a switching
strategy. Both the switching strategy and the
feasibility check utilize   the
errors from switching between the different discretized models, and errors due  different mesh sizes in the discretization.
The plot in Figure~\ref{fig:estimator-error-time} shows a steady decrease of the values of the error {measure} estimators over the course of
the iterations of the adaptive algorithm.

Different third-party optimization
  software packages were  compared with the adaptive algorithm on a part of an
existing real-world district heating network. None of the classical optimization tools
converged to a feasible point while the adaptive solver always
terminated after a finite number of iterations with a locally optimal solution of a model that
has a physical accuracy for which state-of-the-art solvers are not
able to compute a feasible solution at all.
\begin{figure}
  \centering
  \resizebox{0.49\textwidth}{!}{\includegraphics{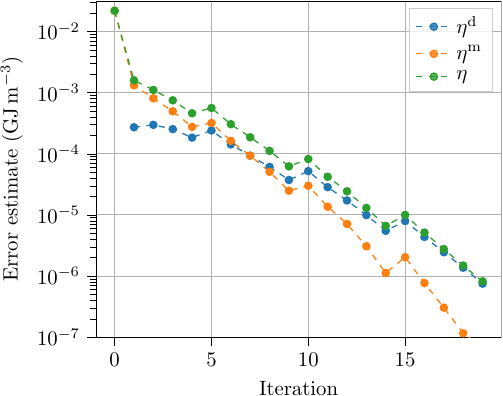}}
  \caption{Error estimator values $\eta$  over the course of the iterations of the adaptive
optimization algorithm using error {measure} estimators. Here $\eta^d$ is the space discretization error and $\eta^m$ the model error estimate.}
  \label{fig:estimator-error-time}
\end{figure}

}

\section*{Conclusions}\label{sec:summary}
In this survey we have discussed automated network based modeling across different physical domains via a system theoretic modeling concept based on differential-algebraic equations, in particular with port-Hamiltonian structure. The approach employs model catalogs to design
space-time-model adaptive methods that can be effectively used in the  co-simulation and optimization of complex dynamical systems, and are essential for designing digital twins.  Modeling, simulation and optimization via port-Hamiltonian descriptor systems is nowadays used in many application domains. We have illustrated the modeling concept for circuit simulation, multi-body dynamics and gas transport as well as district heating networks. However, many challenges remain. These include the development of open source  software packages, the error and perturbation analysis, {large scale optimization methods for real industrial models} and the incorporation of stochastic components.


\begin{thebibliography}{10}

\bibitem{AltGPSS26}
R.~Altmann, I.~C. Garcia, E.~Paakkunainen, P.~Schulze, and S.~Sch{\"o}ps.
\newblock Energy-based modeling for field--circuit coupling.
\newblock {\em Applied Mathematical Modelling}, 155:116688, 2026.

\bibitem{AntBG20}
A.~C. Antoulas, C.~A. Beattie, and S.~Gugercin.
\newblock {\em Interpolatory Methods for Model Reduction}.
\newblock Computational Science \& Engineering. Society for Industrial and
  Applied Mathematics, Philadelphia, PA, USA, 2020.

\bibitem{BanR03}
W.~Bangerth and R.~Rannacher.
\newblock {\em Adaptive finite element methods for differential equations}.
\newblock Springer Science \& Business Media, 2003.

\bibitem{BenCOW17}
P.~Benner, A.~Cohen, M.~Ohlberger, and K.~Willcox.
\newblock {\em Model Reduction and Approximation}.
\newblock Computational Science \& Engineering. Society for Industrial and
  Applied Mathematics, Philadelphia, PA, 2017.

\bibitem{BouVC20}
S.~Bouzid, P.~Viarouge, and J.~Cros.
\newblock Real-time digital twin of a wound rotor induction machine based on
  finite element method.
\newblock {\em Energies}, 13(20):5413, 2020.

\bibitem{BreCP96}
K.~E. Brenan, S.~L. Campbell, and L.~R. Petzold.
\newblock {\em Numerical solution of initial-value problems in
  differential-algebraic equations}.
\newblock Society for Industrial and Applied Mathematics, Philadelphia, PA,
  1996.

\bibitem{Cam87a}
S.~L. Campbell.
\newblock A general form for solvable linear time varying singular systems of
  differential equations.
\newblock {\em {SIAM} J. Math. Anal.}, 18:1101--1115, 1987.

\bibitem{CheMS22}
K.~Cherifi, V.~Mehrmann, and P.~Schulze.
\newblock Simulations in a digital twin of an electrical machine.
\newblock Technical report, 2022.
\newblock http://arxiv.org/abs/2207.02171.

\bibitem{CheSMGL26}
K.~Cherifi, P.~Schulze, V.~Mehrmann, L.~Go{\ss}lau, and P.~L{\"u}nnemann.
\newblock Hierarchical modeling for an industrial implementation of a digital
  twin for electrical drives.
\newblock In K.~Cherifi and I.V. Gosea, editors, {\em Physics-Based and
  Data-Driven Modeling for Digital Twins, ICIAM 2023}. Springer Nature
  Singapore Pte Ltd, 2026.
\newblock http://arxiv.org/abs/2207.02171.

\bibitem{ChuM25b}
D.~Chu and V.~Mehrmann.
\newblock Asymptotic stability and strict passivity of port-{H}amiltonian
  descriptor systems via state feedback.
\newblock {\em Systems and Control Letters}, 202:1006116, 2025.

\bibitem{ChuM25a}
D.~Chu and V.~Mehrmann.
\newblock Stabilization of linear port-{H}amiltonian descriptor systems via
  output feedback.
\newblock {\em {SIAM} J. Matrix Anal. Appl.}, 46:1280--1300, 2025.

\bibitem{DaeMRS24}
H.~D{\"a}nschel, V.~Mehrmann, M.~Roland, and M.~Schmidt.
\newblock Adaptive nonlinear optimization of district heating networks based on
  model and discretization catalogs.
\newblock {\em SEMA Journal}, 81:81--112, 2024.

\bibitem{DomschkePhD}
P.~Domschke.
\newblock {\em Adjoint-Based Control of Model and Discretization Errors for Gas
  Transport in Networked Pipelines}.
\newblock PhD thesis, TU Darmstadt, 2011.

\bibitem{DomDSLM18}
P.~Domschke, A.~Dua, J.~J. Stolwijk, J.~Lang, and V.~Mehrmann.
\newblock Adaptive refinement strategies for the simulation of gas flow in
  networks using a model hierarchy.
\newblock {\em Electron. Trans. Numer. Anal.}, 48:97--113, 2018.
\newblock arxiv: http://arxiv.org/abs/1701.09031.

\bibitem{DomHLMMT21}
P.~Domschke, B.~Hiller, J.~Lang, V.~Mehrmann, R.~Morandin, and C.~Tischendorf.
\newblock Gas network modeling: an overview.
\newblock Preprint, Technische Universit{\"a}t Darmstadt, 2021.
\newblock Available from \emph{https://opus4.kobv.de/opus4-trr154}.

\bibitem{Domschke2011b}
P.~Domschke, O.~Kolb, and J.~Lang.
\newblock Adjoint-based control of model and discretisation errors for gas and
  water supply networks.
\newblock In X.~Yang and S.~Koziel, editors, {\em Computational Optimization
  and Applications in Engineering and Industry}, pages 1--17. Springer, Berlin
  Heidelberg, 2011.

\bibitem{EicF98}
E.~Eich-Soellner and C.~F\"uhrer.
\newblock {\em Numerical Methods in Multibody Dynamics}.
\newblock Vieweg+Teubner Verlag, Wiesbaden, 1998.

\bibitem{DTreport24}
K.~E.~Willcox et~al.
\newblock Foundational research gaps and future directions for digital twins.
\newblock Technical report, Washington, DC: The National Academies Press, 2024.

\bibitem{FalK21}
G.~Falekas and A.~Karlis.
\newblock Digital twin in electrical machine control and predictive
  maintenance: state-of-the-art and future prospects.
\newblock {\em Energies}, 14(18):5933, 2021.

\bibitem{Fre11}
R.~W. Freund.
\newblock The {SPRIM} algorithm for structure-preserving order reduction of
  general {RLC} circuits.
\newblock In P.~Benner, M.~Hinze, and E.~J.~W. ter Maten, editors, {\em Model
  reduction for circuit simulation}, pages 25--52. Springer-Verlag, Dordrecht,
  2011.

\bibitem{GomTBLV18}
C.~Gomes, C.~Thule, D.~Broman, P.~G. Larsen, and H.~Vangheluwe.
\newblock Co-simulation: a survey.
\newblock {\em ACM Computing Surveys (CSUR)}, 51(3):1--33, 2018.

\bibitem{GunBJR21}
M.~G{\"u}nther, A.~Bartel, B.~Jacob, and T.~Reis.
\newblock Dynamic iteration schemes and port-{H}amiltonian formulation in
  coupled differential-algebraic equation circuit simulation.
\newblock {\em Int J. Circ. Theor. Appl.}, 49(2):430--452, 2021.

\bibitem{Hau24}
S.-A. Hauschild.
\newblock {\em Structure-Preserving Numerical Approximations for a
  Port-Hamiltonian Formulation of the Non-Isothermal Euler Equations}.
\newblock PhD thesis, Univ. Trier, 2024.

\bibitem{HauMMMMRS19}
S.-A. Hauschild, N.~Marheineke, V.~Mehrmann, J.~Mohring, A.~Moses Badlyan,
  M.~Rein, and M.~Schmidt.
\newblock Port-hamiltonian modeling of disctrict heating networks.
\newblock {\em DAE Forum, Progress in Differential-Algebraic Equations II},
  2020.
\newblock 333--355.

\bibitem{HauMMMMRS20}
S.-A. Hauschild, N.~Marheineke, V.~Mehrmann, J.~Mohring, A.~Moses~Badlyan,
  M.~Rein, and M.~Schmidt.
\newblock Port-{H}amiltonian modeling of district heating networks.
\newblock In {\em Progress in Differential-Algebraic Equations II}, pages
  333--355, Paderborn, Germany, 2020.

\bibitem{HauM21}
Sarah-Alexa Hauschild and Nicole Marheineke.
\newblock Extended group finite element method for a port-hamiltonian
  formulation of the non-isothermal euler equations.
\newblock {\em PAMM}, 21(1):e202100032, 2021.

\bibitem{JacZ12}
B.~Jacob and H.~Zwart.
\newblock {\em Linear port-{H}amiltonian systems on infinite-dimensional
  spaces}.
\newblock Operator Theory: Advances and Applications. Birkh{\"a}user, Basel,
  2012.

\bibitem{KruMS21}
R.~Krug, V.~Mehrmann, and M.~Schmidt.
\newblock Nonlinear optimization of district heating networks.
\newblock {\em Optimization and Engineering}, 22(2):783--819, 2021.

\bibitem{KunM98}
P.~{Kunkel} and V.~{Mehrmann}.
\newblock Regular solutions of nonlinear differential-algebraic equations and
  their numerical determination.
\newblock {\em Numer. Math.}, 79(4):581--600, 1998.

\bibitem{KunM01}
P.~{Kunkel} and V.~{Mehrmann}.
\newblock Analysis of over- and underdetermined nonlinear
  differential-algebraic systems with application to nonlinear control
  problems.
\newblock {\em Math. Control Signals Systems}, 14(3):233--256, 2001.

\bibitem{KunM18}
P.~{Kunkel} and V.~{Mehrmann}.
\newblock Regular solutions of {DAE} hybrid systems and regularization
  techniques.
\newblock {\em {BIT} Numer. Math.}, 58:1049--1077, 2018.

\bibitem{KunM24}
P.~Kunkel and V.~Mehrmann.
\newblock {\em Differential-Algebraic Equations. Analysis and Numerical
  Solution}.
\newblock EMS Press, Berlin, Germany, 2nd edition, 2024.

\bibitem{LamMT13}
R.~Lamour, R.~M{\"a}rz, and C.~Tischendorf.
\newblock {\em Differential-algebraic equations: a projector based analysis}.
\newblock Differential-Algebraic Equations Forum. Springer-Verlag, Berlin,
  Heidelberg, 2013.

\bibitem{MarLHLPTU26}
A.~Martin, F.~Liers, F.~Hante, J.~Lang, M.~E. Pfetsch, C.~Tischendorf, and
  S.~Ulbrich, editors.
\newblock {\em Mathematical Modelling, Simulation and Optimization using the
  Example of Gas Networks}, volume 174 of {\em International Series of
  Numerical Mathematics}.
\newblock Birkh{\"a}user, Cham, 2026.

\bibitem{MehMW18}
C.~Mehl, V.~Mehrmann, and M.~Wojtylak.
\newblock Linear algebra properties of dissipative {H}amiltonian descriptor
  systems.
\newblock {\em {SIAM} J. Matrix Anal. Appl.}, 39(3):1489--1519, 2018.

\bibitem{Meh15}
V.~{Mehrmann}.
\newblock Index concepts for differential-algebraic equations.
\newblock In {\em Encyclopedia of Applied and Computational Mathematics}, pages
  676--681. Springer-Verlag, Berlin, Heidelberg, 2015.

\bibitem{MehM19}
V.~Mehrmann and R.~Morandin.
\newblock Structure-preserving discretization for port-{H}amiltonian descriptor
  systems.
\newblock In {\em 58th IEEE Conference on Decision and Control (CDC), Nice,
  France}, pages 6863--6868, 2019.

\bibitem{MehSS18}
V.~Mehrmann, M.~Schmidt, and J.J. Stolwijk.
\newblock Model and discretization error adaptivity within stationary gas
  transport optimization.
\newblock {\em Vietnam J. Mathematics}, 46:779--801, 2018.

\bibitem{MehU23}
V.~Mehrmann and B.~Unger.
\newblock Control of port-{H}amiltonian differential-algebraic systems and
  applications.
\newblock {\em Acta Numerica}, To appear, 2023.

\bibitem{MehS23}
V.~Mehrmann and A.J. van~der Schaft.
\newblock Differential-algebraic systems with dissipative hamiltonian
  structure.
\newblock {\em Math. Control Signals Systems}, 2023.

\bibitem{Mor24}
R.~Morandin.
\newblock {\em Modeling and numerics of port-{H}amiltonian descriptor systems}.
\newblock Dissertation, Technische Universit\"at Berlin, 2024.

\bibitem{NedPS22}
N.~Nedialkov, J.~D. Pryce, and L.~Scholz.
\newblock An energy-based, always index $\leq 1$ and structurally amenable
  electrical circuit model.
\newblock {\em SIAM Journal on Scientific Computing}, 44(4):B1122--B1147, 2022.

\bibitem{NocSV09}
R.H. Nochetto, K.G. Siebert, and A.~Veeser.
\newblock Theory of adaptive finite element methods: an introduction.
\newblock In {\em Multiscale, nonlinear and adaptive approximation: Dedicated
  to Wolfgang Dahmen on the occasion of his 60th birthday}, pages 409--542.
  Springer, 2009.

\bibitem{NusT16}
T.~Nussbaumer and S.~Thalmann.
\newblock Influence of system design on heat distribution costs in district
  heating.
\newblock {\em Energy}, 101:496--505, 2016.

\bibitem{Poe22}
B.~P{\"o}chtrager.
\newblock {\em Coupling Multiphysical Systems in Automotive Simulation
  Software/submitted by Dipl.-Ing. Bernhard P{\"o}chtrager}.
\newblock PhD thesis, Institut f\"ur Industriemathematik, Johannes Kepler
  Universit\"at, Linz, Austria, 2022.

\bibitem{PolS10}
R.~V. Polyuga and A.~{van der} Schaft.
\newblock Structure preserving model reduction of port-{H}amiltonian systems by
  moment matching at infinity.
\newblock {\em Automatica J. IFAC}, 46(4):665--672, 2010.

\bibitem{RabR00}
P.~C. Rabier and W.~C. Rheinboldt.
\newblock {\em Nonholonomic motion of rigid mechanical systems from a {DAE}
  viewpoint}.
\newblock Society for Industrial and Applied Mathematics, Philadelphia, PA,
  USA, 2000.

\bibitem{RabR96a}
P.~J. Rabier and W.~C. Rheinboldt.
\newblock Classical and generalized solutions of time-dependent linear
  differential-algebraic equations.
\newblock {\em Linear Algebra Appl.}, 245:259--293, 1996.

\bibitem{RabR96b}
P.~J. Rabier and W.~C. Rheinboldt.
\newblock Time-dependent linear {DAE}s with discontinuous inputs.
\newblock {\em Linear Algebra Appl.}, 247:1--29, 1996.

\bibitem{RasCSS20}
R.~Rashad, F.~Califano, A.~J. van~der Schaft, and S.~Stramigioli.
\newblock {Twenty years of distributed port-{H}amiltonian systems: a literature
  review}.
\newblock {\em {IMA} J. Math. Control I.}, pages 1--23, 2020.

\bibitem{RasO20}
A.~Rasheed, O.~San, and T.~Kvamsdal.
\newblock Digital twin: Values, challenges and enablers from a modeling
  perspective.
\newblock {\em IEEE Access}, 8:21980--22012, 2020.

\bibitem{Ria08}
R.~Riaza.
\newblock {\em Differential-Algebraic Systems. Analytical Aspects and Circuit
  {A}pplications}.
\newblock World Scientific Publishing Co. Pte. Ltd., Hackensack, NJ., 2008.

\bibitem{SchABHJKKOPS17}
M.~Schmidt, D.~A{\ss}mann, R.~Burlacu, J.~Humpola, I.~Joormann, N.~Kanelakis,
  T.~Koch, D.~Oucherif, M.~E. Pfetsch, L.~Schewe, et~al.
\newblock Gaslib—a library of gas network instances.
\newblock {\em Data}, 2(4):40, 2017.

\bibitem{SerMH19}
A.~Serhani, D.~Matignon, and G.~Haine.
\newblock A partitioned finite element method for the structure-preserving
  discretization of damped infinite-dimensional port-{H}amiltonian systems with
  boundary control.
\newblock In F.~Nielsen and F.~Barbaresco, editors, {\em Geometric Science of
  Information}, pages 549--558. Springer, Cham, 2019.

\bibitem{ShaCE22}
V.~Shashkov, I.~Cortes~Garcia, and H.~Egger.
\newblock {MONA}—a magnetic oriented nodal analysis for electric circuits.
\newblock {\em International Journal of Circuit Theory and Applications},
  50(9):2997--3012, 2022.

\bibitem{Sim13}
B.~Simeon.
\newblock {\em Computational flexible multibody dynamics. A
  Differential-Algebraic Approach}.
\newblock Differential-Algebraic Equations Forum. Springer-Verlag, Berlin,
  Heidelberg, 2013.

\bibitem{StoM18}
J.~J. Stolwijk and V.~Mehrmann.
\newblock Error analysis and model adaptivity for flows in gas networks.
\newblock {\em Analele Stiintifice Univ. Ovidius Constanta. Seria Matematica},
  26(2):231--266, 2018.

\bibitem{Tre13}
S.~Trenn.
\newblock Solution concepts for linear {DAE}s: {A} survey.
\newblock In A.~Ilchmann and T.~Reis, editors, {\em Surveys in
  Differential-Algebraic Equations I}, Differential-Algebraic Equations Forum,
  pages 137--172. Springer-Verlag, Berlin, Heidelberg, 2013.

\bibitem{Sch13}
A.~van~der Schaft.
\newblock Port-{H}amiltonian differential-algebraic systems.
\newblock In A.~Ilchmann and T.~Reis, editors, {\em Surveys in
  Differential-Algebraic Equations I}, Differential-Algebraic Equations Forum,
  pages 173--226. Springer-Verlag, Berlin, Heidelberg, 2013.

\bibitem{SchJ14}
A.~{van der S}chaft and D.~Jeltsema.
\newblock Port-{H}amiltonian systems theory: {A}n introductory overview.
\newblock {\em Foundations and Trends in Systems and Control}, 1(2-3):173--378,
  2014.

\bibitem{ZwaM24}
H.~Zwart and V.~Mehrmann.
\newblock Abstract dissipative hamiltonian differential-algebraic equations are
  everywhere.
\newblock {\em DAE Panel}, 2, 2024.

\end{thebibliography}

\end{document}